\documentclass{amsart}
\title[On the vanishing of some mock theta functions]
{On the vanishing of some mock theta functions at odd roots of unity}
\usepackage{amssymb,amsmath,amsthm,epsfig,graphics,latexsym,hyperref}
\theoremstyle{definition}
\newtheorem{definition}{Definition}
\theoremstyle{plain}
\newtheorem{lemma}      {Lemma}

\newtheorem{theorem}    {Theorem}

\newtheorem{corollary}  {Corollary}
\newtheorem{conjecture} {Conjecture}
\newtheorem{open} {Open Problem}
\theoremstyle{remark}

\newcommand{\fr}{\frac}

 \newmuskip\pFqskip
\mathchardef\pFcomma=\mathcode`, % keep a copy of the comma

\begin{document}
  \author[M. El Bachraoui]{Mohamed El Bachraoui}
  \address{Dept. Math. Sci,
 United Arab Emirates University, PO Box 15551, Al-Ain, UAE}
 \email{melbachraoui@uaeu.ac.ae}
% \address{American University of Sharjah, Department of Mathematics and Statistics,
%American University of Sharjah, P.O. Box 26666, Sharjah, UAE}
%\email{jgriffin@aus.edu}
%
\begin{abstract}
We consider the problem of whether or not certain
mock theta functions vanish at the roots of unity with an odd order.
We prove for any such function $f(q)$ that there exists a constant $C>0$ such
that for any odd integer $n>C$ the function $f(q)$ does vanish at the primitive $n$-th roots of
unity. This leads us to conjecture that $f(q)$ does not vanish at the primitive
$n$-th roots of unity for any odd positive integer $n$.
\end{abstract}
\date{\textit{\today}}
\maketitle
\section{Introduction}
Throughout $q$ is a complex number with $|q|<1$, $n$ is a nonnegative integer, and $p$ is an odd prime number. A complex number $\alpha$ is an $n$-th root of unity
if $\alpha^n = 1$ and it is called a primitive $n$-th root of unity if $n$ is the smallest nonnegative integer such that $\alpha^n=1$.
If the order of $\alpha$ is odd we will sometimes say that $\alpha$ is
an odd root of unity.
The asymptotic behaviour at the roots of unity is an important feature in the theory of mock theta functions.
Indeed, for a function $f(q)$ defined by $q$-series to be a mock theta function one of Ramanujan's
conditions states:

($\ast$)\quad
 for every root of unity $\zeta$, there is a theta function $\theta_{\zeta}(q)$ such that the
\par
\qquad difference $f(q)-\theta_{\zeta}(q) = O(1)$ as $q\to \zeta$ radially.

\noindent
The first examples of mock theta functions were given by Ramanujan~\cite{Ramanujan} and there is
nowadays an intensive literature dealing with these functions, relations among them, identities relating
them to other $q$-series, and their asymptotic behaviour at root of unity according to the
condition~($\ast$). For an account of the mock theta functions see for
example~\cite{Andrews-1986, Andrews-Hickerson, Watson-1936, Watson-1937}. Folsom-Ono-Rhoades~\cite{Fol-Ono-Rho-2013} went further with the condition~($\ast$) as
they also gave explicit formulas for $O(1)$ as a linear combination of the root of unity and its powers.
For other references focusing more on the asymptotic behaviour of mock theta functions at the roots of unity including explicit formulas for $O(1)$ in terms of these roots,
we refer to~\cite{BKLMR} and~\cite{Fol-Ono-Rho}.

Our current motivation is of a quite different flavour as we will consider the problem of whether or not
certain mock theta functions have zeros at the odd roots of unity.
By way of example, let $\phi(q)$ and $\sigma(q)$ be the sixth order mock theta functions of
Ramanujan~\cite{Ramanujan} given by
\[
\phi(q) = \sum_{k\geq 0} \frac{(-1)^k q^{k^2} (q;q^2)_k}{(-q;q)_{2k}}\quad \text{and\quad }
\sigma(q) = \sum_{k\geq 0} \frac{q^{\binom{k+2}{2}}(-q;q)_k}{(q;q^2)_{k+1}}
\]
where we have used the standard notation from the theory of
basic hypergeometric series~\cite{Gasper-Rahman}:
\[
(a;q)_0  = 1,\  (a;q)_n = \prod_{j=0}^{n-1} (1-aq^j),\  \text{and\ }
(a;q)_{\infty}  = \prod_{j=0}^{\infty} (1-aq^j).
\]
For convenience we let
\[
(a_1,\ldots, a_k; q)_n = \prod_{j=1}^k (a_j;q)_n \ \text{and\ }
(a_1,\ldots, a_k; q)_{\infty} = \prod_{j=1}^k (a_j;q)_{\infty}.
\]
Rather than $\sigma(q)$ we will consider $\sigma(-q)$ since the latter function is defined
and it terminates at the odd roots of unity.
Ramanujan~\cite{Ramanujan} stated the following formula
which was first proved by
Andrews-Hickerson~\cite[(0.19$)_R$]{Andrews-Hickerson}
\begin{equation*}\label{combination-1}
\phi(q^2) + 2\sigma(-q) = (q;q^2)_{\infty}^2 (q^3;q^6)_{\infty}^2 (q^6;q^6)_{\infty}.
\end{equation*}
Then it is clear from the right hand-side of this identity that the function $\phi(q^2) + 2\sigma(-q)$
vanishes at the $n$-th roots of unity for any odd $n$. So, it is natural to ask the question
whether or not the individual
terms $\phi(q)$ and $\sigma(-q)$ vanish at these $n$-th roots of unity.
Besides, a very important remark is that our functions when evaluated at an odd primitive root of unity
$\zeta$ become finite sums
of the form $f(\zeta)=\sum_{j} c_{\zeta^{j}} \zeta^j$ with rational coefficients $c_{\zeta^j}$.
This is in connection with the theme of vanishing sums of roots of unity as presented in the paper by Conway-Jones~\cite{Conway-Jones}.

Conway-Jones~\cite{Conway-Jones} considered certain equations involving trigonometric functions with rational arguments which they called trigonometric Diophantine equations.
The authors reduced solving these equations to the problem of finding all linear combinations with rational coefficients
among the roots of unity and they were able to classify all such linear combinations of less than ten terms.
In particular, they considered linear relations of the form
\(
\sum_{i=0}^k a_i \zeta^{n_i} = 0
\)
where $a_i$ are rational numbers and $\zeta$ is some root of unity of order $Q$ and showed that for every $C>1$ one has
\[
\log Q \leq C\sqrt{k\log k} + O(1).
\]
As was indicated in~\cite{Conway-Jones}, vanishing sums of roots of unity have been also looked at earlier by R\'{e}dei~\cite{Redei}.
Lenstra~\cite{Lenstra} studied the asymptotic behaviour of the coefficients
occurring in linear combinations with integer coefficients of the roots of unity. For instance, the author proved that
if $\sum_{i=1}^k a_i \zeta_i =0$ where $\zeta_i$ are roots of unity and $a_i$ are integers such that $\gcd(a_1,\ldots,a_k)=1$
and if no proper subset of $\{\zeta_1,\ldots,\zeta_k\}$ is linearly dependent over $\mathbb{Q}$, then
\[
|a_i| \leq 2^{1-k} k^{k/2}.
\]
Zannier~\cite{Zannier}
went further and allowed complex coefficients as he investigated linear relations
\(
\sum_{i=0}^k a_i \zeta^{n_i} = 0
\)
where $a_i$ are complex numbers and $\zeta$ is some root of unity of order $Q$. The author as a main
result obtained bounds for the order $Q$. Lam-Leung~\cite{Lam-Leung} for a given positive integer $m$
classified the positive integers $n$ for which there exist $m$-th roots of unity
$\zeta_1,\ldots, \zeta_n$ such that $\sum_{i=1}^n \zeta_i = 0$.

Our main objective in this work is to prove that for sufficiently large positive odd integer $n$, the functions
$\phi(q)$ and $\sigma(-q)$ along with a variety of other mock theta functions do not vanish at the $n$-th
roots of unity. This leads us to conjecture that these functions do not vanish at the $n$-th roots
of unity for any positive odd integer $n$.

Furthermore, one can also consider for odd $n$ sums
\begin{equation}\label{vanishing-sum}
\sum\{ f(\alpha_j):\ \alpha_j\ \text{is an $n$-th root of unity} \}
\end{equation}
and ask the question for which odd positive integers $n$
the sum~\eqref{vanishing-sum} is equal to zero.
It is of course tempting in the cases at hand to try to come up with closed formulas for $\phi(\zeta)$ and
$\sigma(-\zeta)$ which yield answers to this question.
However, letting $\zeta$ be a $p$-th root of unity for the small prime values $p=3, 5, 7, 11$,
it is easy to verify that:
\begin{equation}\label{small-p}
%\begin{itemize}
\begin{split}
\text{if $p=3$, then\ } & \phi(\zeta) = -2\zeta \ \text{and\ } \sigma(-\zeta)=\zeta^2, \\
\text{if $p=5$, then\ } & \phi(\zeta) = \zeta^4 \ \text{and\ } \sigma(-\zeta)= -\fr{\zeta^3}{2}, \\
\text{if $p=7$, then\ } & \phi(\zeta) = 2-\zeta \ \text{and\ } \sigma(-\zeta)= -1+\fr{\zeta^2}{2}, \\
\text{ if $p=11$, then\ } & \phi(\zeta) = \zeta-\zeta^2-\zeta^4+3\zeta^8  \\
& \text{and\ } \sigma(-\zeta)= -\fr{\zeta^2}{2} + \fr{\zeta^4}{2} -\fr{3\zeta^5}{2} + \fr{\zeta^8}{2},
\end{split}
\end{equation}
which suggests that it would not be easy to find such closed formulas.
In addition, observe from~\eqref{small-p} that if $p=5$ then the sum~\eqref{vanishing-sum} vanishes for
both $f(q)=\phi(q)$ and $f(q)=\sigma(-q)$ since
\[
\begin{split}
\sum\{ \phi(\alpha_j):\ \alpha_j\ \text{is a $5$-th root of unity} \}
&= \phi(1) + \phi(\zeta) + \phi(\zeta^2) + \phi(\zeta^3) + \phi(\zeta^4) \\
&= 1 + \zeta^4 + \zeta^3 + \zeta^2 + \zeta = 0
\end{split}
\]
and
\[
\begin{split}
\sum\{ \sigma(-\alpha_j):\ \alpha_j\ \text{is a $5$-th root of unity} \}
&= \sigma(-1) + \sigma(-\zeta) + \sigma(-\zeta^2) + \sigma(-\zeta^3) + \sigma(-\zeta^4) \\
& = -\fr{1}{2} - \fr{\zeta^3}{2} -\fr{\zeta}{2} - \fr{\zeta^4}{2} - \fr{\zeta^2}{2} = 0
\end{split}
\]
where $\zeta$ is any primitive $5$-th root of unity. However, for the cases $p=3, 7, 11$
the sum~\eqref{vanishing-sum} for these two functions does not vanish.
We note that we do not know the answer to this question and it is not our purpose to discuss it in the present work.
%Note that $f(1)$ is also a term in the sum~\eqref{vanishing-sum}.

In this paper we will pay attention to the vanishing problem at the odd roots of unity for the following
sixth order mock theta functions of Ramanujan~\cite{Ramanujan}
\begin{equation}\label{sixth}
\begin{split}
\phi(q) &= \sum_{k\geq 0} \frac{(-1)^k q^{k^2} (q;q^2)_k}{(-q;q)_{2k}}, \quad
\psi(q) = \sum_{k\geq 0} \frac{(-1)^k q^{(k+1)^2} (q;q^2)_k}{(-q;q)_{2k+1}} ,\\
\lambda(q) &= \sum_{k\geq 0} \frac{(-1)^k q^{k} (q;q^2)_k}{(-q;q)_{k}} ,\quad
\mu(q) = \sum_{k\geq 0} \frac{(-1)^k(q;q^2)_k}{(-q;q)_{k}}, \\
\rho(q) &= \sum_{k\geq 0} \frac{q^{\binom{k+1}{2}}(-q;q)_k}{(q;q^2)_{k+1}} ,\quad
\sigma(q) = \sum_{k\geq 0} \frac{q^{\binom{k+2}{2}}(-q;q)_k}{(q;q^2)_{k+1}},
\end{split}
\end{equation}
the following eighth order mock theta functions of Gordon-McIntosh~\cite{Gordon-McIntosh}
\begin{equation}\label{eighth}
\begin{split}
S_0(q) &= \sum_{k\geq 0}\fr{q^{k^2} (-q;q^2)_k}{(-q^2;q^2)_{k}}, \quad
S_1(q) = \sum_{k\geq 0}\fr{q^{k(k+2)} (-q;q^2)_k}{(-q^2;q^2)_{k}}, \\
U_0(q) &= \sum_{k\geq 0}\fr{q^{k^2} (-q;q^2)_k}{(-q^4;q^4)_{k}}, \quad
U_1(q) = \sum_{k\geq 0}\fr{q^{(k+1)^2} (-q;q^2)_k}{(-q^2;q^4)_{k}} ,
%V_0(q) &= -1 + 2 \sum_{k\geq 1}\fr{q^{k^2} (-q;q^2)_k}{(q;q^2)_{k}} \\
%V_1(q) &=\sum_{k\geq 1}\fr{q^{(k+1)^2} (-q;q^2)_k}{(q;q^2)_{k}}.
\end{split}
\end{equation}
the following fifth order mock theta functions of Ramanujan~\cite{Ramanujan-1927}
\begin{equation}\label{fifth}
\begin{split}
\phi_0(q) &= \sum_{k\geq 0} q^{k^2} (-q;q^2)_k\ \text{and\ }
\phi_1(q) = \sum_{k\geq 0} q^{(k+1)^2} (-q;q^2)_k,
\end{split}
\end{equation}
and finally the following second-order mock theta function of Ramanujan~\cite{Ramanujan}
\begin{equation}\label{second}
u(q)= \sum_{k\geq 0} \fr{(-1)^k q^{k^2} (q;q^2)_k}{(-q^2;q^2)_{k}^2}.
\end{equation}
We restricted ourselves to the mock theta functions in~\eqref{sixth},~\eqref{eighth},~\eqref{fifth}, and~\eqref{second}
because they satisfy the following properties required by our argument of proof. Firstly,
either the series representing $f(q)$ or $f(-q)$ is both defined and terminating at the odd roots of unity. Observe that the terminating condition is guaranteed by the
presence of $(q;q^2)_k$ and $(q;q)_k$ in the numerator of the $k$-th term of our series. So, for instance, our argument does not work for Ramanujan's
third order mock theta function
\[
f(q) = \sum_{n=0}^{\infty}\frac{q^{n^2}}{(-q;q)_n^2}
\]
as $f(q)$ does not terminate at the odd roots of unity and $f(-q)$ has singularities at these roots. A similar restriction applies for the fifth and the tenth order
mock theta functions
\[
\sum_{n=0}^{\infty}\frac{q^{2n^2}}{(q;q^2)_n}\ \text{and\ }
\sum_{n=0}^{\infty}\frac{q^{n(n+1)/2}}{(q;q^2)_n}.
\]
Our third condition is that the function $f(q)$ should be bounded at the odd roots of unity.
Related to this, we recall that Andrews-Hickerson in~\cite[Theorem 5.0]{Andrews-Hickerson}
to confirm the asymptotic condition~($\ast$),
needed to establish the boundedness of the two functions $\phi(q)$ and $\psi(q)$
at the odd roots of unity and Gordon-McIntosh~\cite[p. 332]{Gordon-McIntosh}
needed to prove the boundedness of the functions $S_0(q)$ and
$S_1(q)$ at such roots.
It turns out that $U_0(q)$ is also bounded at the odd roots of unity by virtue of
the following formula of Gordon-McIntosh~\cite{Gordon-McIntosh}:
\begin{equation}\label{U0S0S1}
U_0(q)= S_0(q^2) + q S_1 (q^2).
\end{equation}
On the other hand, as we are not aware whether it is well-known that the functions
$\phi_0(-q)$, $\phi_1(-q)$, and $u(q)$  are bounded at the odd roots of unity, we record
this as a lemma which we will prove later in Section~\ref{sec:proof-lem}.
\begin{lemma}\label{lem-1}
The functions $\phi_0(-q)$, $\phi_1(-q)$, and $u(q)$
are bounded at the primitive odd roots of unity.
\end{lemma}
The following notation is justified.
%As the boundedness property of the functions $\phi(q)$, $\psi(q)$, $S_0(q)$, $S_1(q)$, and $U_0(q)$
%is crucial to our work, we record the following notations for further reference.
\begin{definition}
For any odd primitive root $\zeta$ let
\[
M_{1,\zeta} = \sup\{\phi(\zeta) \},
\
M_{2,\zeta} = \sup\{\psi(\zeta) \},
\]
\[
M_{3,\zeta} = \sup\{S_0(\zeta) \},
\
M_{4,\zeta} = \sup\{S_1(\zeta) \},
\
M_{5,\zeta} = \sup\{U_0(\zeta) \},
\]
\[
M_{6,\zeta} = \sup\{\phi_0(\zeta) \},
\
M_{7,\zeta} = \sup\{\phi_1(\zeta) \},
\ \text{and\ }
M_{8,\zeta}=\sup\{u(\zeta) \}.
\]
Furthermore,  let for $j=1,\ldots,8$
\[
M_j= \inf\{ M_{j,\zeta}:\ \zeta\ \text{is a primitive root of unity of odd order} \} .
%M_2= \inf\{ M_{2,\zeta}:\ \zeta\ \text{is a primitive root of unity} \}.
\]
\end{definition}
The rest of the paper is organized as follows. In Section~\ref{sec:results} we collect our main theorems
and their corollaries along with proofs for these corollaries. Sections~\ref{sec:proof-1}-\ref{sec:proof-5}
 are devoted to proofs of the main theorems and in Section~\ref{sec:proof-lem} we give the proof of
 Lemma~\ref{lem-1}.
Finally, in Section~\ref{sec:concluding} we shall give some comments and state our conjectures and
open problems.
\section{Main results}\label{sec:results}
In this section we collect our main theorems along with their corollaries.
\begin{theorem}\label{thm phi-mu}
If $n$ is an odd positive integer such that $n >M_1^2$, then the functions $\phi(q)$, $\mu(q)$, and
$\sigma(-q)$ do not vanish at the primitive $n$-th roots of unity.
\end{theorem}
\begin{corollary}\label{cor phi-mu}
If $n$ is an odd positive integer such that $n >M_1^2$, then the functions $\phi'(q)$, $\mu'(q)$, and
$\sigma'(-q)$ defined by
\[
\begin{split}
\phi'(q) &= \sum_{k\geq 0} \fr{q^k (q;q^2)_k}{(-q;q)_{2k}} \\
\mu'(q) &=  1+ \sum_{k\geq 1} \fr{q^{-\binom{k-1}{2}} (q;q^2)_k}{(-q;q)_k} \\
\sigma'(-q) &=  \sum_{k\geq 0} \fr{(-1)^k (q;q)_k}{(-q;q^2)_{k+1}}
\end{split}
\]
do not vanish at the primitive $n$-th roots of unity.
\end{corollary}
\begin{proof}
As
\[
\begin{split}
(q^{-1};q^{-2})_k &= (1-\fr{1}{q}) (1-\fr{1}{q^3})\cdots (1-\fr{1}{q^{2k-1}}) \\
&= (-1)^k q^{-k^2} (q;q^2)_k
\end{split}
\]
and
\[
\begin{split}
(-q^{-1};q^{-1})_{2k} &= (1+\fr{1}{q}) (1+\fr{1}{q^2})\cdots (1+\fr{1}{q^{2k}}) \\
&= q^{-k(2k+1)} (-q;q)_{2k}
\end{split}
\]
we get
\[
\fr{(-1)^k q^{-k^2} (q^{-1};q^{-2})_k}{(-q^{-1};q^{-1})_{2k}}
= \fr{q^k (q;q^2)_k}{(-q;q)_{2k}}
\]
and therefore we deduce that $\phi'(q) = \phi(q^{-1})$. So, assuming that $\phi'(\zeta)=0$
for a primitive $n$-th root of unity means that $\phi(\zeta^{-1})=0$ which contradicts
Theorem~\ref{thm phi-mu} since $\zeta^{-1}$ is also a primitive $n$-th root of unity. The same argument applies
to $\mu'(q)= \mu(q^{-1})$ and $\sigma'(-q)=\sigma(q^{-1})$.
\end{proof}
\begin{theorem}\label{thm psi-lambda}
If $n$ is an odd positive integer such that $n >4 M_2^2$, then the functions $\psi(q)$, $\lambda(q)$,
and $\rho(-q)$ do not vanish at the primitive $n$-th roots of unity.
\end{theorem}
\begin{corollary}\label{cor psi-lambda}
If $n$ is an odd positive integer such that $n >4M_2^2$, then the functions $\psi'(q)$, $\lambda'(q)$, and
$\rho'(-q)$ defined by
\[
\begin{split}
\psi'(q) &= \sum_{k\geq 0} \fr{q^k (q;q^2)_k}{(-q;q)_{2k+1}} \\
\lambda'(q) &=  \sum_{k\geq 0} \fr{q^{-\binom{k}{2}} (q;q^2)_k}{(-q;q)_k} \\
\rho'(-q) &=  \sum_{k\geq 0} \fr{(-1)^k q^{k+1}(q;q)_k}{(-q;q^2)_{k+1}}
\end{split}
\]
do not vanish at the primitive $n$-th roots of unity.
\end{corollary}
\begin{proof}
Proofs follow easily from Theorem~\ref{thm psi-lambda} and the relations
\[
\psi'(q) = \psi(q^{-1}),\
\lambda'(q) = \lambda(q^{-1}),\ \text{and\ }
\rho'(-q) = \rho(q^{-1}).
\]
\end{proof}
\begin{theorem}\label{thm S0S1}
Let $n$ be an odd positive integer.

(a)\ If $n > M_3^2$, then the function $S_0(-q)$
does not vanish at the primitive $n$-th roots of unity.

(b)\ If $n > M_4^2$, then the function $S_1(-q)$
does not vanish at the primitive $n$-th roots of unity.

(c)\ If $n > M_5^2$, then the functions $U_0(-q)$
and $U_1(-q)$ do not vanish at the primitive $n$-th roots of unity.
\end{theorem}
\begin{corollary}\label{cor S0S1}
Let $n$ be an odd positive integer and let
\[
\begin{split}
S_0'(q) & = \sum_{k\geq 0} \fr{(-1)^k q^{k(k+1)} (-q;q^2)_k}{(-q^2;q^2)_{k}} \\
S_1'(q) &= \sum_{k\geq 0} \fr{(-1)^k q^{k(k-1)} (-q;q^2)_k}{(-q^2;q^2)_{k}} \\
U_0'(q) &= \sum_{k\geq 0} \fr{(-1)^k q^{2k(k+1)} (-q;q^2)_k}{(-q^4;q^4)_{k}} \\
U_1'(q) &= \sum_{k\geq 0} \fr{(-1)^{k+1} q^{2k(k+1)+1} (-q;q^2)_k}{(-q^2;q^4)_{k}}.
\end{split}
\]

(a)\ If $n >M_3^2$, then the function $S_0'(-q)$  does not vanish at the primitive $n$-th roots of unity.

(b)\ If $n >M_4^2$, then the function $S_1'(-q)$  does not vanish at the primitive $n$-th roots of unity.

(c)\ If $n >M_5^2$, then the functions $U_0'(-q)$  and $U_1'(-q)$ do not vanish at the primitive
$n$-th roots of unity.
\end{corollary}
\begin{proof}
Part~(a) follows by Theorem~\ref{thm S0S1}(a) combined with the identity
$S_0'(q) = S_0(q^{-1})$, part~(b) by Theorem~\ref{thm S0S1}(b) combined with the identity
$S_1'(q) = S_1(q^{-1})$, and part~(c) is a consequence of Theorem~\ref{thm S0S1}(c) and
the relations $U_0'(q) = U_0(q^{-1})$ and $U_1'(q) = U_1(q^{-1})$.
\end{proof}
\begin{theorem}\label{thm phi0-1}
Let $n$ be an odd positive integer.

(a)\ If $n > M_6^2$, then the function $\phi_0(-q)$
does not vanish at the primitive $n$-th roots of unity.

(b)\ If $n > M_7^2$, then the function $\phi_1(-q)$
does not vanish at the primitive $n$-th roots of unity.
\end{theorem}
\begin{theorem}\label{thm u}
Let $n$ be an odd positive integer.
If $n > M_8^2$, then the function $u(q)$
does not vanish at the primitive $n$-th roots of unity.
\end{theorem}
%
%\section{Sums vanishing at roots of unity}
%
\section{Proof of Theorem~\ref{thm phi-mu}}\label{sec:proof-1}
Let $n$ be an odd positive integer.
It is easily verified that $\phi(q)$ does not vanish at the $n$-th roots of unity for $n=1, 3, 5$.
Now suppose that $n>\sup\{5,M_1^2\}$ and let $\zeta$ be a primitive $n$-th root of unity. Then
\begin{equation}\label{help-1-1}
\begin{split}
\phi(\zeta)
&= \sum_{k=0}^{\fr{n-1}{2}} (-1)^k \fr{ \zeta^{k^2}(\zeta;\zeta^2)_k}{(-\zeta;\zeta)_{2k}} \\
&= 1- \fr{\zeta(1-\zeta)}{(1+\zeta)(1+\zeta^2)}
+ \fr{\zeta^4 (1-\zeta)(1-\zeta^3)}{(1+\zeta)(1+\zeta^2)(1+\zeta^3)(1+\zeta^4)} \\
&\quad
+\cdots + (-1)^{\fr{n-1}{2}} \fr{\zeta^{(\fr{n-1}{2})^2} (\zeta;\zeta^2)_{\fr{n-1}{2}}}{(-\zeta;\zeta)_{n-1}}.
\end{split}
\end{equation}
Since for any primitive $n$-th root of unity $\alpha$ we have the following basic relation
\begin{equation}\label{basic-1}
(-\alpha;\alpha)_{n-1}=(1+\alpha)(1+\alpha^2)\cdots (1+\alpha^{n-1}) =1,
\end{equation}
the identity~\eqref{help-1-1} becomes
\begin{equation}\label{help-1-2}
\begin{split}
\phi(\zeta)
&= 1-\zeta(1-\zeta)(1+\zeta^3)\cdots (1+\zeta^{n-1})
+ \zeta^4(1-\zeta)(1-\zeta^3)(1+\zeta^5)\cdots (1+\zeta^{n-1}) \\
&\quad + \cdots+ (-1)^{\fr{n-1}{2}} \zeta^{(\fr{n-1}{2})^2} (\zeta;\zeta^2)_{\fr{n-1}{2}}.
\end{split}
\end{equation}
From
\[
\zeta(1-\zeta^{n-1}) = \zeta -1,\
\zeta^3(1-\zeta^{n-3}) = \zeta^3 -1,\ \ldots,\
\zeta^{n-2}(1-\zeta^{2}) = \zeta^{n-2} -1
\]
we get
\[
(1-\zeta^2)(1-\zeta^4)\cdots (1-\zeta^{n-1}) = (-1)^{-\fr{n-1}{2}}  \zeta^{-(\fr{n-1}{2})^2}
(1-\zeta)(1-\zeta^3)\cdots (1-\zeta^{n-2}).
\]
This combined with the elementary fact
\[
(\zeta;\zeta)_{n-1} = (1-\zeta)(1-\zeta^2)\cdots (1-\zeta^{n-1}) = n
\]
yields
\[
\Big( (1-\zeta)(1-\zeta^3)\cdots (1-\zeta^{n-2}) \Big)^2 (-1)^{-\fr{n-1}{2}} \zeta^{-(\fr{n-1}{2})^2} = n.
\]
Then
\begin{equation}\label{basic-2}
(-1)^{\fr{n-1}{2}} \zeta^{(\fr{n-1}{2})^2} (\zeta;\zeta^2)_{\fr{n-1}{2}}
=  (-1)^{\fr{n-1}{4}} \zeta^{\fr{(n-1)^2}{8}} \sqrt{n}
\end{equation}
and thus
\begin{equation}\label{key}
|(-1)^{\fr{n-1}{2}} \zeta^{(\fr{n-1}{2})^2} (\zeta;\zeta^2)_{\fr{n-1}{2}} | = \sqrt{n}.
\end{equation}
Using~\eqref{help-1-2} and~\eqref{key} we achieve
\[
\begin{split}
| \phi(\zeta)| &=  \Big|\sum_{k=0}^{\fr{n-3}{2}}(-1)^k \fr{\zeta^{k^2} (\zeta;\zeta^2)_{k}}{(-\zeta;\zeta)_{2k}}
+ (-1)^{\fr{n-1}{2}} \zeta^{(\fr{n-1}{2})^2} (\zeta;\zeta^2)_{\fr{n-1}{2}} \Big| \\
& \geq \left| \Big|\sum_{k=0}^{\fr{n-3}{2}}(-1)^k \fr{\zeta^{k^2} (\zeta;\zeta^2)_{k}}{(-\zeta;\zeta)_{2k}} \Big|
-|(-1)^{\fr{n-1}{2}} \zeta^{(\fr{n-1}{2})^2} (\zeta;\zeta^2)_{\fr{n-1}{2}} | \right| \\
& = \sqrt{n} - \Big|\sum_{k=0}^{\fr{n-3}{2}}(-1)^k \fr{\zeta^{k^2} (\zeta;\zeta^2)_{k}}
{(-\zeta;\zeta)_{2k}}\Big| \\
& \geq \sqrt{n} - M_1 \\
& > 0.
\end{split}
\]
Consequently, $\phi(q)$ does not vanish at the primitive $n$-th roots of unity for any odd $n>M_1^2$. As for $\mu(q)$, we need the following identity of
Ramanujan~\cite{Ramanujan} which
was first confirmed by Andrews-Hickerson~\cite[(0.20$)_R$]{Andrews-Hickerson}
\begin{equation*}
2 \phi(q^2) - 2\mu(-q) = (-q;q^2)_{\infty}^2 (-q^3;q^6)_{\infty}^2 (q^6;q^6)_{\infty}.
\end{equation*}
Replacing $q$ by $-q$ in the previous relation we get
\begin{equation}\label{help phi-mu}
2 \phi(q^2) - 2\mu(q) = (q;q^2)_{\infty}^2 (q^3;q^6)_{\infty}^2 (q^6;q^6)_{\infty}.
\end{equation}
Observe that the difference $\phi(q^2)- \mu(q)$ vanishes at the $n$-th
roots of unity since the right hand-side of~\eqref{help phi-mu}
clearly does. Now suppose for a contradiction that
$\mu(\zeta)=0$ for some primitive $n$-th roots of unity. Then so does $\phi(\zeta^2)$ by the previous observation.
But this contradicts the first statement of the theorem since $\zeta^2$ is also a primitive $n$-th root
of unity.
Similarly, for the statement regarding $\sigma(-q)$ we appeal the following identity
of Ramanujan~\cite{Ramanujan} which was first proved by
Andrews-Hickerson~\cite[(0.19$)_R$]{Andrews-Hickerson}
\begin{equation*}
\phi(q^2) + 2\sigma(q) = (-q;q^2)_{\infty}^2 (-q^3;q^6)_{\infty}^2 (q^6;q^6)_{\infty}.
\end{equation*}
Upon replacing $q$ by $-q$ we find
\begin{equation}\label{help phi-sigma}
\phi(q^2) + 2\sigma(-q) = (q;q^2)_{\infty}^2 (q^3;q^6)_{\infty}^2 (q^6;q^6)_{\infty}
\end{equation}
from which it follows that $\phi(q^2) + 2\sigma(-q)$ vanishes at the $n$-th roots of unity.
Thus assuming that $\sigma(-\zeta)=0$ at some primitive $n$-th root of unity implies that
$\phi(\zeta^2)=0$ which is impossible by the first statement of the theorem.
\section{Proof of Theorem~\ref{thm psi-lambda}}\label{sec:proof-2}
We proceed in the same way as in the proof of Theorem~\ref{thm phi-mu}.
Let $n$ be an odd positive integer. It is easy to check that
$\psi(q)$ does not vanish at the $n$-th roots of unity for $n=1, 3, 5$.
Let $n>\sup\{5,4M_2^2\}$ and let
 $\zeta$ be a primitive $n$-th root of unity. We have with the help of~\eqref{basic-1} and~\eqref{basic-2}
\begin{equation}\label{help-2-1}
\begin{split}
\psi(\zeta)
&= \sum_{k=0}^{\fr{n-1}{2}} (-1)^k \fr{ \zeta^{(k+1)^2}(\zeta;\zeta^2)_k}{(-\zeta;\zeta)_{2k+1}} \\
&= \fr{1}{1+\zeta}- \fr{\zeta^4(1-\zeta)}{(1+\zeta)(1+\zeta^2)(1+\zeta^3)} \\
&\quad + \fr{\zeta^9 (1-\zeta)(1-\zeta^3)}{(1+\zeta)(1+\zeta^2)(1+\zeta^3)(1+\zeta^4)(1+\zeta^5)} \\
& \quad
+\cdots+ (-1)^{\fr{n-1}{2}} \fr{\zeta^{(\fr{n+1}{2})^2} (\zeta;\zeta^2)_{\fr{n-1}{2}}}{(-\zeta;\zeta)_{n-1} (1+1)}\\
&=(1+\zeta^2)\cdots(1+\zeta^{n-1}) -\zeta^4(1-\zeta)(1+\zeta^4)(1+\zeta^5)\cdots (1+\zeta^{n-1}) \\
&\quad
+\cdots + (-1)^{\fr{n-1}{2}}\fr{\zeta^{(\fr{n+1}{2})^2}}{2} (\zeta;\zeta^2)_{\fr{n-1}{2}} \\
&=(1+\zeta^2)\cdots(1+\zeta^{n-1}) -\zeta^4(1-\zeta)(1+\zeta^4)(1+\zeta^5)\cdots (1+\zeta^{n-1}) \\
&\quad
+\cdots +  \fr{\sqrt{n}}{2} (-1)^{\fr{n-1}{4}}\zeta^{\fr{(n+1)^2 + 4n}{8}}.
\end{split}
\end{equation}
Now combine~\eqref{help-2-1} with~\eqref{key} to obtain
\[
\begin{split}
| \psi(\zeta)| &=  |\sum_{k=0}^{\fr{n-3}{2}}(-1)^k \fr{\zeta^{(k+1)^2} (\zeta;\zeta^2)_{k}}{(-\zeta;\zeta)_{2k+1}}
+ (-1)^{\fr{n-1}{2}} \zeta^{(\fr{n+1}{2})^2} (\zeta;\zeta^2)_{\fr{n+1}{2}} | \\
& \geq \fr{\sqrt{n}}{2} -
|\sum_{k=0}^{\fr{n-3}{2}}(-1)^k \fr{\zeta^{(k+1)^2} (\zeta;\zeta^2)_{k}}{(-\zeta;\zeta)_{2k+1}}| \\
& \geq \fr{\sqrt{n}}{2} - M_2 \\
& > 0.
\end{split}
\]
This shows that $\psi(q)$ does not vanish at the primitive $n$-th roots of unity.
Regarding the function
$\lambda(q)$ we first make an appeal to the following formula which was stated
by Ramanujan~\cite{Ramanujan} and proved by
Andrews-Hickerson~\cite[(0.21)$_R$]{Andrews-Hickerson}
\begin{equation*}
2 q^{-1}\psi(q^2) + \lambda(-q) = (-q;q^2)_{\infty}^2 (-q,-q^5,q^6;q^6)_{\infty} .
\end{equation*}
Replacing $q$ by $-q$ we get
\begin{equation}\label{help psi-lambda}
-2 q^{-1}\psi(q^2) + \lambda(q) = (q;q^2)_{\infty}^2 (q,q^5,q^6;q^6)_{\infty}
\end{equation}
from which we deduce that $-2 q^{-1}\psi(q^2) + \lambda(q)$ vanishes at the primitive roots of unity
of odd order. So, assuming  that $\lambda(\zeta)=0$ for some primitive $n$-th root of unity $\zeta$ implies that
$\zeta^{-1}\psi(\zeta^2)=0$. A contradiction since $\zeta^2$ is also a primitive $n$-th root of unity.
As for the function $\rho(q)$, we combine similar ideas with the following identity
of Ramanujan~\cite{Ramanujan} which was confirmed in~\cite[(0.18)$_R$]{Andrews-Hickerson}
\begin{equation}\label{help psi-rho}
q^{-1}\psi(q^2) + \rho(q) = (-q;q^2)_{\infty}^2 (-q,-q^5,q^6;q^6)_{\infty} .
\end{equation}
\section{Proof of Theorem~\ref{thm S0S1}}\label{sec:proof-3}
Let $n$ be an odd positive integer.
As the result is evident for $n=1,3,5$ we assume that $n > \sup\{5, M_3^2\}$. Let
$\zeta$ be a primitive $n$-th root of unity. From
\[
S_0(-q) = \sum_{k\geq 0} \fr{(-1)^k q^{k^2} (q;q^2)_k}{(-q^2;q^2)_k},
\]
we clearly see that the sum on the right hand-side terminates at $\zeta$. Since $\zeta^2$ is also a
primitive $n$-th root of unity, we obtain with the help of~\eqref{basic-1}
\[
\begin{split}
S_0(-\zeta) &= \sum_{k=0}^{\fr{n-3}{2}} \fr{(-1)^k \zeta^{k^2} (\zeta;\zeta^2)_k}{(-\zeta^2;\zeta^2)_k}
+ \fr{(-1)^{\fr{n-1}{2}}\zeta^{(\fr{n-1}{2})^2} (\zeta;\zeta^2)_{\fr{n-1}{2}}}{(-\zeta^2;\zeta^2)_{\fr{n-1}{2}}} \\
&= 1 - \zeta(1-\zeta)(1+\zeta^4)(1+\zeta^6)\cdots (1+\zeta^{2n-2}) \\
&+ \zeta^4(1-\zeta)(1-\zeta^3)(1+\zeta^6)(1+\zeta^8)\cdots (1+\zeta^{2n-2}) \\
&+\cdots + (-1)^{\fr{n-3}{2}}\zeta^{(\fr{n-3}{2})^2} (\zeta;\zeta^2)_{\fr{n-3}{2}}
(1+\zeta^{n-1})(1+\zeta^{n+1})\cdots (1+\zeta^{2n-2}) \\
& + \fr{(-1)^{\fr{n-1}{2}}\zeta^{(\fr{n-1}{2})^2} (\zeta;\zeta^2)_{\fr{n-1}{2}}}{(-\zeta^2;\zeta^2)_{\fr{n-1}{2}}}.
\end{split}
\]
We now claim that the norm of the $(n-1)/2$-term of $S_0(-\zeta)$ is
\[
\left| \fr{(-1)^{\fr{n-1}{2}}\zeta^{(\fr{n-1}{2})^2} (\zeta;\zeta^2)_{\fr{n-1}{2}}}{(-\zeta^2;\zeta^2)_{\fr{n-1}{2}}}
\right| = \sqrt{n}.
\]
From the formulas
\[
\zeta^2(1+\zeta^{n-2}) = 1+\zeta^2,\ \zeta^4(1+\zeta^{n-4}) = 1+\zeta^4,\ldots,
\zeta^{n-1}(1+\zeta^{2}) = 1+\zeta^{n-1}
\]
we deduce that
\[
(1+\zeta)(1+\zeta^3)\cdots (1+\zeta^{n-2}) =
\zeta^{-\fr{n^2 -1}{4}}(1+\zeta^2)(1+\zeta^4)\cdots (1+\zeta^{n-1})
\]
which combined with~\eqref{basic-1} gives
\[
1= \Big( (1+\zeta^2)(1+\zeta^4)\cdots (1+\zeta^{n-1}) \Big)^2 \zeta^{-\fr{n^2 -1}{4}}
\]
or equivalently,
\begin{equation}\label{key-1}
(-\zeta^2;\zeta^2)_{\fr{n-1}{2}} = \zeta^{\fr{n^2 -1}{8}}.
\end{equation}
Now use the foregoing formula and~\eqref{key} to derive the claim.
Then
\[
\begin{split}
| S_0(-\zeta)| &=  \left| \sum_{k=0}^{\fr{n-3}{2}}(-1)^k \fr{\zeta^{k^2} (\zeta;\zeta^2)_{k}}
{(-\zeta^2;\zeta^2)_{k}}
+ \fr{(-1)^{\fr{n-1}{2}} \zeta^{(\fr{n-1}{2})^2} (\zeta;\zeta^2)_{\fr{n-1}{2}}}
{(-\zeta^2;\zeta^2)_{\fr{n-1}{2}}} \right| \\
& \geq \sqrt{n} - \left| \sum_{k=0}^{\fr{n-3}{2}}(-1)^k \fr{\zeta^{k^2} (\zeta;\zeta^2)_{k}}
{(-\zeta^2;\zeta^2)_{k}} \right| \\
& \geq \sqrt{n} - M_3 \\
& > 0
\end{split}
\]
which means that $S_0(-\zeta) \not= 0$ for any primitive $n$-th root of unity $\zeta$.
This proves part~(a).
Proofs for parts (b) and (c) regarding the functions $S_1(-q)$ and $U_1(-q)$ are omitted as they follow similarly.
As for $U_1(-q)$, we need the following formula of Gordon-McIntosh~\cite[(1.9)]{Gordon-McIntosh}
(with $q$ replaced by $-q$)
\begin{equation}\label{U0U1}
U_0(-q)+2U_1(-q)= (q;q^2)_{\infty}^3 (q^2;q^2)_{\infty} (q^2;q^4)_{\infty}.
\end{equation}
Clearly the sum $U_0(-q)+2U_1(-q)$  vanishes at the $n$-th roots of unity.
Assume for a contradiction that $U_1(-\zeta)=0$ for some primitive $n$-th root of unity. Then obviously
$2 U_1(-\zeta)=0$. Now combining with~\eqref{U0U1} gives $U_0(-\zeta)=0$,
which is absurd.
\section{Proof of Lemma~\ref{lem-1}}\label{sec:proof-lem}
The idea of proof for all the three functions is essentially the same as the proof of~\cite[Theorem 5.0]{Andrews-Hickerson}
but we establish only the statement on the function $u(q)$ as it requires more effort.
Consider the following eighth order mock theta functions from~\cite{Gordon-McIntosh}
\[
\begin{split}
T_0(q) &= \sum_{k\geq 0} \fr{q^{(k+1)(k+2)} (-q^2;q^2)_k}{(-q;q^2)_{k+1}} \\
T_1(q) &= \sum_{k\geq 0} \fr{q^{k(k+1)} (-q^2;q^2)_k}{(-q;q^2)_{k+1}}.
\end{split}
\]
We note that in~\cite[p.332]{Gordon-McIntosh} it has been shown that $T_0(q)$ and $T_1(q)$ are both bounded at the even roots of unity.
We claim that $T_1(q)$ is bounded at the primitive odd roots of unity. Let $\zeta$ be a primitive $N$-th root of unity for an odd positive integer $N$ and write $q=r\zeta$ for
$0\leq r \leq 1$. As the claim is clear for $r=0$, we assume that $0<r \leq 1$. Let
\[
u_n(r) = \fr{q^{n(n+1)} (-q^2;q^2)_n}{(-q;q^2)_{n+1}}.
\]
Then it is easy to see that
\begin{equation}\label{help-lem-1}
|u_{n+N}(r)| = q^{2Nn + N(N+1)} \fr{|(-q^{2n+2};q^2)_N|}{|(q^{2n+3};q^2)_N|} |u_n(r)|.
\end{equation}
By an application of~\cite[Lemma 5.2]{Andrews-Hickerson} to the foregoing identity with $R=r^{2n+2k}$ and $R'=r^{2n+2N}$, we get
\[
\begin{split}
|(-q^{2n+2};q^2)_N| &= \prod_{k=1}^N |1+r^{2n+2k} \zeta^{2n+2k}| \\
&\leq \prod_{k=1}^N r^{k-N} |1+r^{2n+2N} \zeta^{2n+2k}| \\
&= r^{N(1-N)/2} |1+r^{2n+2N} \zeta^{2n+2k}|.
\end{split}
\]
Since $\zeta$ is a primitive $N$-th root of unity and $k$ runs from $1$ to $N$, we have that
$\zeta^{2n+2k}$ runs through the $N$-th roots of unity and therefore $1+r^{2n+2N} \zeta^{2n+2k}$
runs through the solutions of the equation
\[
(x-1)^N - r^{N(2n+2N)}.
\]
Then by Vieta's formula, the product of these solutions is
$(-1)^N (1-r^{N(2n+2N)})$
and thus we obtain
\begin{equation}\label{help-lem-2}
| (-q^{2n+2};q^2)_N| \leq r^{N(1-N)/2} (1-r^{N(2n+2N)}).
\end{equation}
In addition, it has been shown in~\cite[p. 95]{Andrews-Hickerson}
that
\[
| (q^{2n+3};q^2)_N| \geq r^{N(N-1)/2}.
\]
Then a combination of the triangular inequality and the previous inequality yields
\begin{equation}\label{help-lem-3}
| (-q^{2n+3};q^2)_N| \geq | (q^{2n+3};q^2)_N| \geq r^{N(N-1)/2}.
\end{equation}
Now combining~\eqref{help-lem-1},~\eqref{help-lem-2}, and~\eqref{help-lem-3} we find
\[
| u_{n+N}(r) | \leq q^{2Nn + 2N} (1-r^{N(2n+2N)}) |u_n(r)|,
\]
which by virtue of~\cite[Lemma 5.3]{Andrews-Hickerson} yields
\[
\begin{split}
| u_{n+N}(r) | &\leq \fr{N(2n+2N)}{2Nn + 2N + N(2n+2N)} |u_n(r)| \\
&\leq \fr{n+N}{2n+N} |u_n(r)| \\
&\leq \fr{8}{11} |u_n(r)|,
\end{split}
\]
provided $n\geq N$. Then by~\cite[Lemma 5.1]{Andrews-Hickerson}, $\sum_{n\geq 0}|u_n(r)|$ is bounded and thus $\sum_{n\geq 0}u_n(r)$ is bounded too. This confirms the claim.
Similarly one can prove that $T_0(q)$ is bounded too at the odd primitive roots of unity.
Moreover, as by Gordon-McIntosh~\cite{Gordon-McIntosh}
\[
U_1(q) = T_0(q^2) + q T_1(q^2),
\]
we see that $U_1(q)$ is also bounded at the odd roots of unity. Recalling that $U_0(q)$ is bounded at these roots of unity and
combining with the following relation of McIntosh~\cite{McIntosh}
\[
u(q) = U_0(q) - 2U_1(q),
\]
we deduce that $u(q)$ is also bounded at the odd roots of unity. This completes the proof.
\section{Proof of Theorem~\ref{thm phi0-1}}\label{sec:proof-4}
We use the same approach as in Section~\ref{sec:proof-1}. From the identity
\begin{equation*}\label{help-4-1}
\begin{split}
\phi_0(-\zeta)
&= \sum_{k=0}^{\fr{n-1}{2}} (-1)^k \zeta^{k^2}(\zeta;\zeta^2)_k \\
&= 1- \zeta(1-\zeta)
+ \zeta^4 (1-\zeta)(1-\zeta^3) \\
&\quad
+\cdots + (-1)^{\fr{n-1}{2}} \zeta^{(\fr{n-1}{2})^2} (\zeta;\zeta^2)_{\fr{n-1}{2}}
\end{split}
\end{equation*}
and~\eqref{key} we get
\[
\begin{split}
| \phi_0(-\zeta)| &=  \Big|\sum_{k=0}^{\fr{n-3}{2}}(-1)^k \zeta^{k^2} (\zeta;\zeta^2)_{k}
+ (-1)^{\fr{n-1}{2}} \zeta^{(\fr{n-1}{2})^2} (\zeta;\zeta^2)_{\fr{n-1}{2}} \Big| \\
& \geq \left| \Big|\sum_{k=0}^{\fr{n-3}{2}}(-1)^k \zeta^{k^2} (\zeta;\zeta^2)_{k} \Big|
-|(-1)^{\fr{n-1}{2}} \zeta^{(\fr{n-1}{2})^2} (\zeta;\zeta^2)_{\fr{n-1}{2}} | \right| \\
& = \sqrt{n} - \Big|\sum_{k=0}^{\fr{n-3}{2}}(-1)^k \zeta^{k^2} (\zeta;\zeta^2)_{k} \Big| \\
& \geq \sqrt{n} - M_6 \\
& > 0.
\end{split}
\]
This shows that $\phi_0(-q)$ does not vanish at the primitive $n$-th roots of unity for any odd $n>M_6^2$.
The proof for the statement on $\phi_1(-q)$ is omitted as it follows similarly.
\section{Proof of Theorem~\ref{thm u}}\label{sec:proof-5}
Noticing the resemblance between $u(q)$ and $S_0(-q)$ and using the same steps as in
Section~\ref{sec:proof-3}, we get
\[
\begin{split}
u_0(\zeta) &= \sum_{k=0}^{\fr{n-3}{2}} \fr{(-1)^k \zeta^{k^2} (\zeta;\zeta^2)_k}{(-\zeta^2;\zeta^2)_k^2}
+ \fr{(-1)^{\fr{n-1}{2}}\zeta^{(\fr{n-1}{2})^2} (\zeta;\zeta^2)_{\fr{n-1}{2}}}{(-\zeta^2;\zeta^2)_{\fr{n-1}{2}}^2} \\
&= 1 - \zeta(1-\zeta)\Big((1+\zeta^4)(1+\zeta^6)\cdots (1+\zeta^{2n-2})\Big)^2 \\
&+ \zeta^4(1-\zeta)(1-\zeta^3) \Big((1+\zeta^6)(1+\zeta^8)\cdots (1+\zeta^{2n-2})\Big)^2 \\
&+\cdots + (-1)^{\fr{n-3}{2}}\zeta^{(\fr{n-3}{2})^2} (\zeta;\zeta^2)_{\fr{n-3}{2}}
\Big( (1+\zeta^{n-1})(1+\zeta^{n+1})\cdots (1+\zeta^{2n-2}) \Big)^2 \\
& + \fr{(-1)^{\fr{n-1}{2}}\zeta^{(\fr{n-1}{2})^2} (\zeta;\zeta^2)_{\fr{n-1}{2}}}
{\Big( (-\zeta^2;\zeta^2)_{\fr{n-1}{2}}\Big)^2}.
\end{split}
\]
Appealing to~\eqref{key} and~\eqref{key-1} we see that the norm of the $(n-1)/2$-term of $u(\zeta)$ is
\[
\left| \fr{(-1)^{\fr{n-1}{2}}\zeta^{(\fr{n-1}{2})^2} (\zeta;\zeta^2)_{\fr{n-1}{2}}}
{\Big( (-\zeta^2;\zeta^2)_{\fr{n-1}{2}} \Big)^2}
\right| = \sqrt{n}.
\]
This yields
\[
\begin{split}
| u(\zeta)| &=  \left| \sum_{k=0}^{\fr{n-3}{2}}(-1)^k \fr{\zeta^{k^2} (\zeta;\zeta^2)_{k}}
{(-\zeta^2;\zeta^2)_{k}^2}
+ \fr{(-1)^{\fr{n-1}{2}} \zeta^{(\fr{n-1}{2})^2} (\zeta;\zeta^2)_{\fr{n-1}{2}}}
{(-\zeta^2;\zeta^2)_{\fr{n-1}{2}}^2} \right| \\
& \geq \sqrt{n} - \left| \sum_{k=0}^{\fr{n-3}{2}}(-1)^k \fr{\zeta^{k^2} (\zeta;\zeta^2)_{k}}
{(-\zeta^2;\zeta^2)_{k}^2} \right| \\
& \geq \sqrt{n} - M_8 \\
& > 0.
\end{split}
\]
This completes the proof.
\section{Concluding remarks}\label{sec:concluding}
In this section we
let $f(q)$ denote any one of the functions in the statements of Theorems~\ref{thm phi-mu}--\ref{thm u},
along with the statements of the listed corollaries.
Based on our theorems and their corollaries stating that for sufficiently large odd $n$, the function $f(q)$ does not vanish at the $n$-th roots of unity, 
we are led to the following conjecture.
\begin{conjecture}
For any odd positive integer and any primitive $n$-th root of unity $\zeta$ we have
$f(\zeta) \not= 0$.
\end{conjecture}
As we observed earlier in the comments just after~\eqref{small-p}, the sum~\eqref{vanishing-sum} for $f(q)=\phi(q)$
and $f(q)=\sigma(-q)$ vanishes at $p=3$ but not at $p=3,7,11$. We have the following open problems.
\begin{open}
Find prime numbers for which the sum~\eqref{vanishing-sum} is zero.
\end{open}
\begin{open}
Is the number of primes for which the sum~\eqref{vanishing-sum} vanishes finite or infinite?
\end{open}
\bigskip
\noindent
{\bf Acknowledgment.} The author is grateful for the referee for valuable comments and suggestions which improved the quality and the presentation of the paper. The final draft of this work has been completed 
in 2021 while visiting my birthplace Nador, Morocco. I thank NOVO~CLASS~Caf\'{e}~\& Restaurant in Nador and its staff for the great environment, the high speed internet, and the top quality coffee.
\end{document}